\documentclass[final,twocolumn,3p]{elsarticle}

\usepackage[utf8]{inputenc}
\usepackage[english]{babel}	
\usepackage[T1]{fontenc}

\usepackage{amsmath,amssymb,amsfonts,amsthm,graphics,mathtools,url, enumerate} 
\usepackage{subcaption}
\usepackage{float}
\usepackage{tikz-cd}

\numberwithin{equation}{section}

\theoremstyle{plain}
\newtheorem{thm}{Theorem}[section]

\theoremstyle{definition}
\newtheorem{ex}[thm]	{Example}
\newtheorem{df}[thm]{Definition}

\newtheorem{rem}[thm]{Remark}

\newcommand{\x}{\times}

\renewcommand{\phi}{\varphi}

\DeclareMathOperator{\arccosh}{arccosh}

\newcommand{\N}{\mathbb{N}}
\newcommand{\Z}{\mathbb{Z}}
\newcommand{\R}{\mathbb{R}}

\newcommand{\bH}{\mathbb{H}}
\newcommand{\bS}{\mathbb{S}}

\begin{document}
\begin{frontmatter}

\title{Integral Betti signature confirms the hyperbolic geometry of brain, climate, and financial networks}
\journal{}

\author[bol]{Luigi Caputi}
\ead{luigi.caputi@unibo.it}
\author[cas,nudz,cvut]{Anna Pidnebesna}
\ead{pidnebesna@cs.cas.cz}
\author[cas,nudz]{Jaroslav Hlinka}
\ead{hlinka@cs.cas.cz}
	
	\address[bol]{Departmente of Mathematics, University of Bologna, Bologna}
\address[cas]{Institute of Computer Science of the Czech Academy of Sciences, \\ Pod Vod\'{a}renskou v\v{e}\v{z}\'{i} 271/2, 182 07 Prague, Czech Republic}
\address[nudz]{National Institute of Mental Health, Topolov\'{a} 748, 250 67 Klecany, Czech Republic}

\begin{abstract}
This paper extends the possibility 
to examine the underlying curvature of data through the lens of topology by using the \emph{Betti curves}, tools of Persistent Homology, as key topological descriptors, building on the clique topology approach. 
It was previously shown that Betti curves distinguish random from Euclidean geometric matrices - i.e.~distance matrices of points randomly  distributed in a cube with Euclidean distance. 
In line with previous experiments, we consider their low-dimensional approximations named \emph{integral Betti values}, or \emph{signatures} that effectively distinguish not only Euclidean, but also spherical and hyperbolic geometric matrices, both from purely random matrices as well as among themselves.
To prove this, we analyse the behaviour of Betti curves for various geometric matrices -- i.e.~distance matrices of points randomly distributed on manifolds of constant sectional curvature, considering the classical models of curvature 0, 1, -1, given by the  Euclidean space, the sphere, and the hyperbolic space. 
We further investigate the dependence of integral Betti signatures on factors including the sample size and dimension. This is important for assessment of real-world connectivity matrices, as we show that the standard approach to network construction gives rise to (spurious) spherical geometry, with topology dependent on sample dimensions.
Finally, we use the  manifolds of constant curvature as comparison models to infer curvature underlying real-world  datasets coming from neuroscience, finance and climate. Their associated topological features exhibit a hyperbolic character: the integral Betti signatures associated to these datasets sit in between Euclidean and hyperbolic (of small curvature). 
The potential confounding ``hyperbologenic effect'' of intrinsic low-rank modular structures is also evaluated through simulations.
\end{abstract}

\end{frontmatter}

{T}he scientific understanding of phenomena around us depends on devising a formal model of reality, which is iteratively tested and further developed based on comparison with observations, and which, ultimately, provides prediction of further phenomena. Many real-world systems, such as the human brain or the Earth's climate, fall within the category of complex systems, with structure neither random nor fully ordered. While the particular structure of such systems may differ from realization to realization, and may also dynamically change  in time, the key challenge is to understand the organizing principles that gave rise to the structure, that is, to 
build a generative model. 

Most often, real-world systems are of very high dimension, and thus, even after some initial dimensionality reduction, they call for representations in high dimensional spaces. Working with such high dimensional objects, although being possible, is often problematic both due to the computational demands, and to the limitations of reliability of such fine-grained analysis. Thus, additional low dimensionality reduction algorithms are usually employed. 
There is a strong and commonly accepted intuition that many real-world, high-dimensional, datasets have a lower dimensional 
representation;  
this is an assumption, usually referred to as the ``manifold hypothesis'': \emph{data in the form of point clouds in~$\R^n$ are sampled from (or, essentially, close to) a manifold~$M^d$ whose dimension~$d$ is smaller than the ambient dimension~$n$}. The manifold $M$ is called \emph{data manifold},  and the recent need of new methodologies for analysing high dimensional data based on this hypothesis gave birth to a new research field, known as  ``manifold learning''.  Experiments show that the manifold hypothesis is true for many real world datasets, 
and  algorithms to test it have been introduced~\cite{MR3522608}. For instance, 
the space of natural images is surprisingly well approximated by a two-dimensional representation homeomorphic to the ($2$-dimensional) 
Klein bottle~\cite{MR3715451};  showing that the underlying topology of data manifolds can be quite unintuitive.
Estimating and understanding the intrinsic geometrical and topological properties of data manifolds is important not just from a theoretical point of view, but also for extracting features and qualitative information from the data itself; e.g.~in developing faster dimensionality-reduction algorithms. 

 The advantage of the manifold hypothesis  is that it allows us to exploit additional theoretical structures and geometric properties of manifolds. Among others, Riemannian geometry is recently gaining much attention~\cite{doi:10.1073/pnas.2100473118}. The usual linear tools, 
 as Principal Component Analysis and Factor Analysis,  
 work well when the data lies close to a linear, flat, subspace of~$\R^n$. However, such linear methods do not work equally well when the data lies near a more complicated, nonlinear,  manifold, and may fail to recover the intrinsic structure, along with related relevant information. Driven by these considerations, numerous nonlinear manifold learning methods and algorithms have been recently proposed; see, e.g.~\cite{geommmod} and the references therein. The notion of (Riemannian) \emph{curvature}
gives an intrinsic, strong, geometric, description of a (data) manifold and has led to the so-called Riemannian Manifold Learning~\cite{Rmanlearn}. 

We refrain here from giving a technical account of techniques  in Riemannian Manifold Learning. The main driving idea is the belief that data might actually be governed by a non-Euclidean, curved, geometry, rather than the most conventional, Euclidean,  one. Developing non-Euclidean methods is therefore of great importance in testing and dealing with such frameworks. 
It is well-known that  Riemannian manifolds of constant (sectional) curvature, in any dimension~$d$, are classified in three classes: of positive, vanishing, and negative curvature. Models of these spaces are the classical sphere~$\bS^d$, the Euclidean space~$\R^d$, and the hyperbolic space~$\bH^d$, respectively. A geometric intuition of how the local intrinsic geometry changes in these models is given by the behaviour of (geodesic) triangles -- see Figure~\ref{fig:triangles}. We refer to~\cite{MR1744486} for an advanced discussion on more detailed geometric properties. Being  positively curved, as opposed to being negatively curved, means that the (data) manifold ``looks like'' a sphere, as opposed to a hyperbolic space, and this comes along with particular (typically) local structural properties, and, ultimately, also functional consequences for the system at hand. On the other hand, most real-world networks seem to have an underlying hyperbolic geometry, which is believed to be related to the fact that ``networks that conform to hyperbolic geometry are maximally responsive to external and internal perturbations''~\cite{SHARPEE2019101}. Hyperbolic geometry, in fact, is thought to be related to faster information transmission \cite{SHARPEE2019101, PhysRevE.82.036106} and some real-world networks have already shown hyperbolic features; e.g.~financial markets~\cite{hypfinance},  olfactory systems~\cite{doi:10.1126/sciadv.aaq1458}, or brain-to-brain coordination networks~\cite{10.3389/fphy.2018.00007}. Also functional brain networks were found to be best represented in hyperbolic spaces~\cite{Whi2021.03.25.436730}. 

Our analysis  starts from these considerations, and aims to provide  qualitative reasoning and tools justifying  the observation that some complex systems have an intrinsic hyperbolic geometry. Investigations of data manifolds, from a Riemannian geometric perspective, need to  \emph{assume} a priori a certain geometry of the underlying data and, as in the case of the sphere, also a specific topology. This leads to the natural questions: what are the main differences in the topological features of data with intrinsic Euclidean and/or non-Euclidean underlying structure? And can we use this knowledge to (at least qualitatively) infer that certain real-world data have intrinsic non-Euclidean geometry, as it has been reported?

To investigate these questions, there are many possible approaches. In~\cite{doi:10.1073/pnas.2100473118}, for example,  methods from differential geometry were borrowed, and the authors aim to compute the curvature of data from the second fundamental form (after additional preprocessing and neighbourhood selection). In the current paper instead, we borrow tools from a recent  subfield of Algebraic Topology, called  Topological Data Analysis (TDA)~\citep{TDAreview1},  hence, combining geometric and topological methods.  
The goal is achieved by employing  manageable topological invariants  as homology and homotopy groups. The structure of the  data is then qualitatively and quantitatively assessed  in the form of topological  features (e.g.~connected components, voids, tunnels or loops). 
 \emph{Persistent Homology} (PH) is one of the main adopted tools, in TDA, to make assessments precise. PH is a multiscale adaptation of the classical homology theories, and it allows a computation of the \emph{persistent} topological features of a space at all different resolutions, while revealing the most essential ones.  It has provided novel qualitative and quantitative approaches to study complex data (in the form of either point clouds, time series or connectivity networks).
Its direct use, along with its derived features, has already been beneficial  in many fields, like neuroscience~\citep{lee-discriminativePH, ourpaper}, finance~\citep{gidea}, image classification~\citep{vmv.20171272}, to name a few. 

In \cite{Giusti13455}, a new PH-derived tool to \emph{reliably detect signatures of structure and randomness that are invariant under non-linear monotone transformations} was introduced, called 
\emph{Betti curves}. 
These features encode information of the topological and potentially also of the geometrical properties of data, represented by (weighted) graphs, as  functions -- curves -- of the edge densities. Our main interest in Betti curves  comes from the fact that they  distinguish geometric Euclidean networks from the random ones~\cite{Giusti13455}. Furthermore, in dimension~$\leq 2$, 
persistent homology and Betti curves detect the curvature~\cite{Bubenik_2020, otter_curv, hacquard2023euler}
of the underlying manifold, and can  be used for obtaining bounds on the Gromov-Hausdorff distance~\cite{Lim2020}. 
Inspired by these recent advancements, in this work we first explore and investigate the behaviour of Betti curves associated to  random sample points uniformly distributed on the standard Riemannian manifolds (in spaces of arbitrary dimension) of sectional curvature $1,0,-1$, and we compare them with the Betti curves of random graphs. For Betti numbers in dimensions~$0$ and $1$, as usually considered in applications, we show that the Betti curves do distinguish the three different geometries, hence the underlying manifolds. Aiming to use  this gained qualitative information in understanding data manifolds, we proceed with analysing three main datasets, coming from brain, stocks, and climate data. These data usually come in the form of correlation matrices, hence we also analyse correlation networks. We observe a general hyperbologenic character, with brain and climate data derived features in between of those derived from the  Euclidean and the hyperbolic geometry of small curvature; the topological features associated to the stocks data having an even more pronounced tendency  towards hyperbolic behaviour. This observation is consistent with the general belief that these data, and in particular the stock market, have an intrinsic hyperbolic geometry. 
 To gain more insight into how the observed results could have been affected by other factors, such as the modular structure -- cf.~\cite{curto2021betti} -- 
 we include a study of the dependence of Betti curves on the number of modules. We observe a qualitative trend going from a Euclidean-like behaviour, in the case of a single module, to hyperbolic, when considering a few modules, to spherical. 
 
 To conclude, we believe that  the methods investigated in this work are effective  descriptors of the underlying geometric manifolds, complementing alternative approaches to analyse the connections between PH and the underlying curvature on manifolds~\cite{fernandez2023intrinsic, curto2021betti} and foresee that quantitative methods based on these approaches will gradually infiltrate data analytic practice similarly as other tools of TDA already have.

 
 We first  give a description of the datasets employed in this work, and then a brief account of the mathematical methods. For each dataset (in the form of time series), we apply the pipeline described in Section~\ref{sec:Betti} to their associated (Pearson) correlation graphs.
 
\section{Data description}
\label{sec:data}

\subsection{Brain data}
\label{sec:data_fmri}

We analyse the dataset that consists of fMRI recordings of 90 healthy controls. 
Functional MRI data were collected with a 3T MR scanner (Siemens; Magnetom Trio) at Institute of Clinical and Experimental Medicine in Prague ($T2^*$-weighted images with BOLD contrast, voxel size 3$\times$3$\times$3mm$^3$, TR/TE = $2000/30$ms, $400$ time points). 
Initial data preprocessing was done using FSL routines (FMRIB Software Library v5.0, Analysis Group, FMRIB, Oxford, UK) and CONN toolbox (McGovern Institute for Brain Research, MIT, USA). 
See~\cite{Kopal2020, ourpaper} for detailed prepocessing descriptions.
To extract the time series for further analysis, the brain’s spatial domain was divided into 90 non-overlapping regions of interest (ROIs) according to the Automatic Anatomical Labelling (AAL) atlas~\cite{Tzourio-Mazoyer2002}; from each ROI we extract one BOLD time series by averaging the time series of all voxels in the ROI.

\subsection{Climate data}
\label{sec:data_climate}

 We use the daily surface air temperature anomalies data obtained from the NECP/NCAR reanalysis dataset~\cite{Kalnay1996}. In particular, we use the daily air temperature fields at 1000hPa level, spanning the period from 1/1/1948  to 12/31/2012 and sampled at  $2.5^\circ\times 2.5^\circ$ angularly regular Gaussian grid. For a more precise description of the dataset, we refer to~\cite{Hlinka2013}. 
 The resulting time series has a length of $23376$ time points at each of $162$ geodesic grid nodes. We randomly select $90$ out of the $162$ time series for better comparability with the brain data.
 
 \subsection{Stocks data}
We use historical stock prices, downloaded from the Yahoo! Finance service~\cite{YahooFinance2023}, belonging to the New York Stock Exchange 100 (NYSE100) index.
We consider only stocks traded between 4 January 1977 and 6 October 2023. This restriction results in $N = 90$ stocks. For daily data, this leads to a data length of $T = 11791$. 
We used the daily adjusted closing prices. The {\em logarithmic return} is  computed as
the first differences of the log-transformed prices:
$r_i(t) = \log\left[ \frac{p_i(t)}{p_i(t-1)}\right],
$
where $p_i(t)$ is the (adjusted, daily) closing price of stock~$i$ at time~$t$.
For a more detailed analysis of a similar dataset, see~\cite{Hartman2018}.
 
 \section{On Betti curves of symmetric matrices}\label{sec:Betti}
 
The topological pipeline employed in this work was first introduced by Giusti \emph{et al.} in~\cite{Giusti13455}.  For a given symmetric (real-valued) matrix, we consider certain topological descriptors  called \emph{Betti curves} -- cf.~Figure~\ref{fig:Drawing}. In order to concisely represent Betti curves, and for comparisons, we use their associated \emph{Area Under the Curve} (AUC),  and call the corresponding features  \emph{integral Betti signatures.} 

\subsection{Betti curves}\label{sec:BettiA}
Betti curves are topological descriptors  of symmetric matrices.
To each  $N\x N$ symmetric matrix~$M$ with distinct non-zero real-valued entries, we first associate a sequence of graphs
\begin{equation}\label{eq:ordcompl}
\mathrm{ord}(M)\coloneqq G_0\subseteq G_1\subseteq \dots\subseteq G_k \ ,
\end{equation}
called  the \emph{order complex}~$\mathrm{ord}(M)$  of $M$ \cite[Def.~2, SI]{Giusti13455}. The construction of $\mathrm{ord}(M)$ starts with a totally disconnected graph~$G_0$ on~$N$ vertices, and proceeds step by step by adding new edges; addition of edges follows the values of the entries of~$M$, sorted either min to max or max to min. We refer to the supplementary notes for a more precise description of the order complex; note that  in all the applications below, values are sorted from the maximum to the minimum value. To each graph~$G_i$ in~$\mathrm{ord}(M)$, we  associate a simplicial complex~$\widetilde{G_i}$ whose $p$-simplices consist of all the complete subgraphs on $(p+1)$  vertices of $G_i$. Classical topological invariants of simplicial complexes are the  \emph{homology groups} -- see, e.g.~\cite{hatcher} and the supplementary -- along with their ranks, called \emph{Betti numbers}. The number of $p$-holes of 
a simplicial complex $\Sigma$ is related to  the $p$-th  {Betti number}~$\beta_p(\Sigma)$ associated to $\Sigma$. Specifically, the $0$-th Betti number~$\beta_0$ gives the number of connected components, and the $1$-st Betti number~$\beta_1$ the number of independent loops.  

For an order complex $\mathrm{ord}(M)$, the pipeline  yields, for each $i\in \N$, a sequence of $i$-th Betti numbers, called the $i$-th \emph{Betti curves} of~$M$:
  \begin{equation}\label{eq:Betti curves}
  	\beta_i(\widetilde{G_0}), \beta_i(\widetilde{G_1}), \dots, \beta_i(\widetilde{G_k}) \ .
  \end{equation}
  	
  	\subsection{Random, geometric and correlation matrices}

We call \emph{random} any symmetric matrix with identically independent real-valued entries, \emph{geometric} any symmetric matrix which is obtained as the distance matrix of sample points randomly distributed on a manifold, and \emph{correlation matrix} any symmetric matrix obtained as the (Pearson) correlation matrix of given time series. 
When dealing with geometric and correlation matrices, we shall consider 
the case of points sampled from the fundamental metric spaces $\R^n$ (Euclidean space), $\mathbb{S}^n$ (a sphere), and $\bH^n$ (a hyperbolic space). It is  known that these are the geometric models of the simply connected $n$-manifolds of constant scalar curvature~$0$, $1$ and $-1$, respectively. 

\begin{figure}[t]
	\centering
	\begin{subfigure}[b]{0.25\textwidth}
		\centering
		\begin{tikzpicture}[baseline=(current bounding box.center),scale =.4]
			\tikzstyle{point}=[circle,thick,draw=black,fill=black,inner sep=0pt,minimum width=2pt,minimum height=2pt]
			\tikzstyle{arc}=[shorten >= 8pt,shorten <= 8pt,->, thick]
			
			\node (v0) at (0,0) {};
			\draw[fill] (0,0)  circle (.05);
			\node (v1) at (3.0,0) {};
			\draw[fill] (3.0,0)  circle (.05);
			\node (v2) at (1.5,2) {};
			\draw[fill] (1.5,2)  circle (.05);
			
			\draw[thick, red] (v0) -- (v1);
			\draw[thick, red] (v0) -- (v2);
			\draw[thick, red] (v1) -- (v2);
		\end{tikzpicture}
		\caption{Triangle in $\R^n$.}\label{fig:YA}
	\end{subfigure}
	\hfill
	\begin{subfigure}[b]{0.15\textwidth}
		\centering
		\begin{tikzpicture}[baseline=(current bounding box.center),scale =.35]
			\tikzstyle{point}=[circle,thick,draw=black,fill=black,inner sep=0pt,minimum width=2pt,minimum height=2pt]
			\tikzstyle{arc}=[shorten >= 8pt,shorten <= 8pt,->, thick]
			
			\node (v0) at (0,0) {};
			\draw[fill] (0,0)  circle (.05);
			\node (v1) at (3.0,0) {};
			\draw[fill] (3.0,0)  circle (.05);
			\node (v2) at (1.5,2) {};
			\draw[fill] (1.5,2)  circle (.05);

			\draw[thick, red] (v0) to[bend right] (v1);
			\draw[thick, red] (v0) to[bend left] (v2);
			\draw[thick, red] (v1) to[bend right] (v2);
		\end{tikzpicture}
		\caption{Triangle in $\mathbb{S}^n$.}\label{fig:YB}
	\end{subfigure}
	\hfill
	\begin{subfigure}[b]{0.15\textwidth}
		\centering
		\begin{tikzpicture}[baseline=(current bounding box.center),scale =.45]
			\tikzstyle{point}=[circle,thick,draw=black,fill=black,inner sep=0pt,minimum width=2pt,minimum height=2pt]
			\tikzstyle{arc}=[shorten >= 8pt,shorten <= 8pt,->, thick]
			
			\node (v0) at (0,0) {};
			\draw[fill] (0,0)  circle (.05);
			\node (v1) at (3.0,0) {};
			\draw[fill] (3.0,0)  circle (.05);
			\node (v2) at (1.5,2) {};
			\draw[fill] (1.5,2)  circle (.05);

			\draw[thick, red] (v0) to[bend left] (v1);
			\draw[thick, red] (v0) to[bend right] (v2);
			\draw[thick, red] (v1) to[bend left] (v2);
		\end{tikzpicture}
		\caption{Triangle in $\mathbb{H}^n$.}\label{fig:YC}
	\end{subfigure}
	\caption{Standard geodesic triangles in  Euclidean, spherical and hyperbolic space.}
	\label{fig:triangles}
\end{figure}
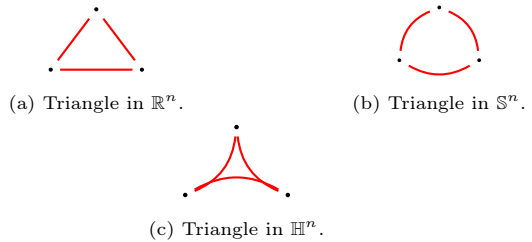

Generally speaking, the geometric properties of a manifold, endowed with different underlying metrics, can be very different. First and foremost,  geodesic triangles in the three fundamental models, have  different behaviours: if  triangles in the Euclidean space are ``straight'', they get ``fatter'' on a sphere and ``thinner'' in a hyperbolic space -- cf.~Figure~\ref{fig:triangles}. An analytical description of the metrics in these geometries is given in the Supplementary. In practical examples, and in our results, we consider sample points randomly distributed in the unit cube, on the unit sphere or in
the hyperbolic space and then compute the Betti curves described in Section~\ref{sec:BettiA}.

\subsection{Dimensional reduction by AUC} 
In order to compare the Betti curves of different datasets, we consider two intrinsic geometric measures: the Area Under the Curve (AUC) of the Betti curves~$\beta_0$ and $\beta_1$. In fact, as shown in Figure~\ref{fig:Drawing}, the qualitative behaviours of the Betti numbers in dimension $i\geq 1$ are all similar, whereas in dimension $0$ they are not. The AUC has been chosen for taking into account the trend of the number of $k$-cycles in a given dimension~$k$, up to affine transformations (as the AUC is invariant under shifts). We will call them 
\emph{integral Betti signatures}. 
Therefore, for each dataset, at different sample sizes, we compute the 
integral Betti signatures, 
obtaining  $3$-dimensional representations -- see Figure~\ref{fig:VDL2d_REALvsSim}.


\section{Results}

The behaviour of the Betti curves associated to \emph{random} graphs has been thoroughly investigated, and it is nowadays theoretically understood -- see, e.g.~\cite{KAHLE20091658}. A more geometric source of (symmetric, random) matrices is given by  \emph{geometric matrices}: instead of considering matrices with random entries, one first considers uniformly distributed random points lying on a geometric manifold (e.g.~the hypercube~$I^n$ or the sphere); then, one  considers the distance matrix~$A$ obtained using the underlying geodesic distances between points~$\{x_1,\dots,x_N\}$, i.e.~setting $A_{i,j}\coloneqq d(x_i,x_j)$. Betti curves associated to sample points in Euclidean spaces  have been investigated by Giusti \emph{et al.}~\citep{Giusti13455}. They show that the Betti curves can distinguish random matrices from the geometric (Euclidean) ones. 

\subsection{Geometric analysis}

\begin{figure}[h!]
\begin{center}
        \includegraphics[width=1.\columnwidth]{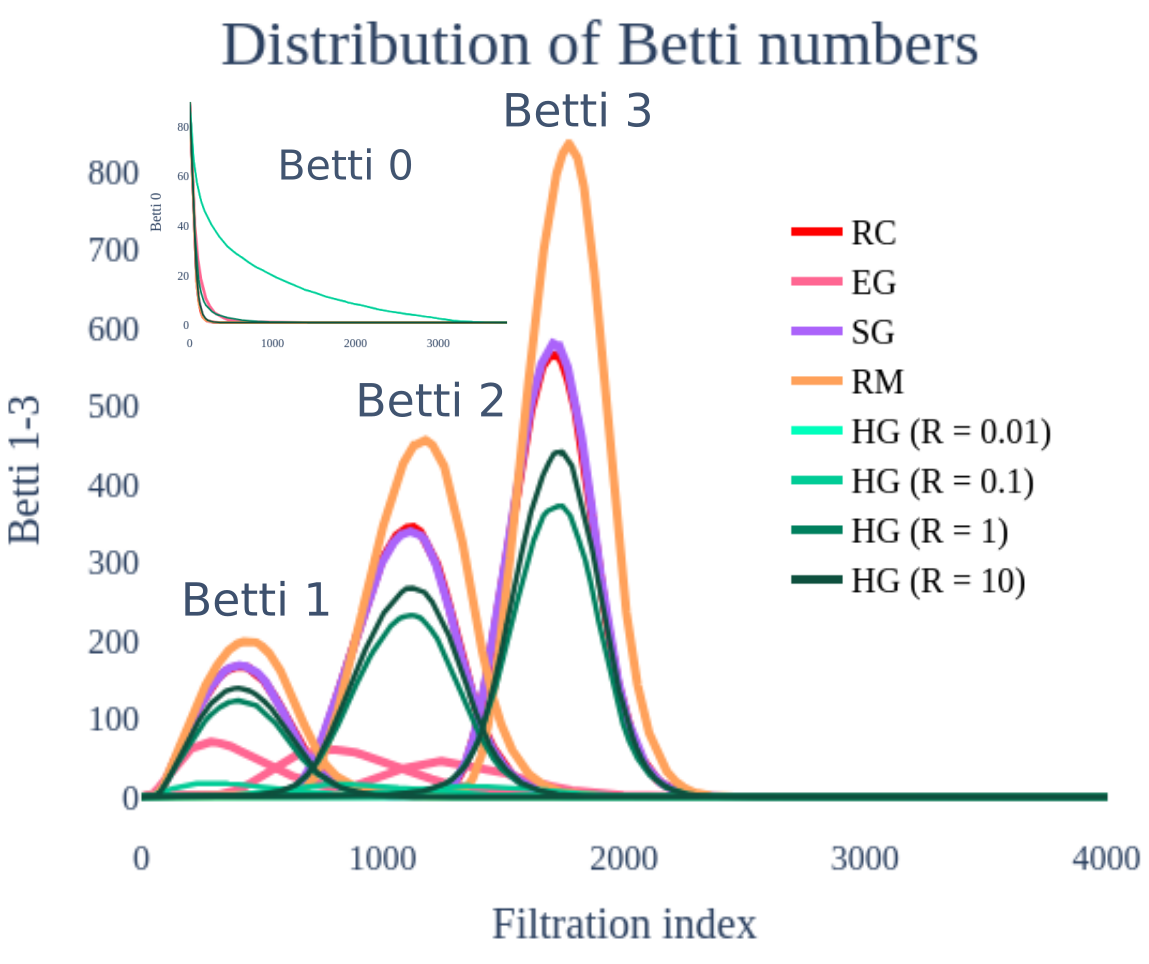}
	\caption{Comparison of Betti numbers distributions of simulated random matrices (RM), Euclidean geometry (EG), spherical geometry (SG), hyperbolic geometry (HG), and random correlation (RC) matrices. The sample size/dimensionality is 400. The Betti curves corresponding to the random correlation are almost not visible on the graph as it underlays the spherical geometry curve. }
	\label{fig:Drawing}
\end{center}
\end{figure}

We start our analysis by extending the results of~\cite{Giusti13455} to  spherical geometry (i.e.~points randomly distributed on the sphere~$\bS^n$, endowed with the spherical  distance) or hyperbolic geometry (i.e.~points randomly  distributed on the hyperbolic space~$\bH^n$ with the hyperbolic distance).  In order to randomly distribute points on such manifolds, we use uniform distributions; a uniform distribution of points in the hypercube $I^n$ in the Euclidean case, and a uniform distribution on the sphere in the spherical case. For the hyperbolic case, we use the Poincar\'e disc model, and the approximation of the distribution -- see, e.g.~\cite{AlanisLobato2016DistanceDB}  -- at given radii. Similar choices were also investigated in \cite{doi:10.1126/sciadv.aaq1458}.  In Figure~\ref{fig:Drawing}, we plot the obtained Betti curves in dimension $0,1,2,3$ associated to random distance matrices (Euclidean: EG, spherical: SG, hyperbolic at different radii: HG), random matrices RM and random correlation matrices~RC (i.e.~sample correlation matrices of Gaussian white noise samples).

\subsubsection*{Dependence on dimension}

As apparent from Figure~\ref{fig:Drawing},  Betti curves distinguish not just random from Euclidean geometric matrices, but they are good descriptors for  other geometries as well. Figure~\ref{fig:Drawing} represents in fact the average Betti curves (over $100$ random iterations) in dimension $0, 1, 2$, and~$3$, of $90\times 90$ random matrices, and geometric matrices of $90$ sample points randomly distributed in the mentioned geometric manifolds of dimension~$400$. In the case of the hyperbolic space, we considered radii $R=0.01, 0.1,1,10$ -- see~\cite{AlanisLobato2016DistanceDB} and the discussion below about the role of $R$. We complemented it with the Betti curves of random correlation matrices with $90$ time series of length~$400$ -- the typical length of fMRI data. 
Under these choices, the higher-dimensional ($\geq 1$) Betti curves associated to random matrices show the highest peaks, followed by the Betti curves of points randomly distributed on the sphere first, in the Euclidean space then, and in the hyperbolic space last. 

Note that the described qualitative behaviour holds across a range of  dimensions of the manifolds (which corresponds to the time series length or sample size) and the number of points sampled, although for small dimensions the order can be reversed. In \cite{Giusti13455}, it was observed that 
 correlation matrices of finite samples of $N$ independent uniformly distributed random variables display the same characteristic Betti curves as random symmetric $N\times N$ matrices (note that for correlating spike trains, the mean of a cross-correlogram was used by~\cite{Giusti13455}). 
However, in our results (for time series  of length $400$) illustrated in Figure~\ref{fig:Drawing}, 
the characteristics of random correlation (RC) matrices do not correspond to the random matrices (RM).
We suggest that the observation presented in \cite{Giusti13455} holds only for sufficiently high dimensions (they used 10000 samples), an observation which might be relevant for the interpretation of TDA analysis of small datasets common in many disciplines.

In Figure~\ref{fig:VDL2d_REALvsSim}, we show  the 3-dimensional plot  in the coordinates (B0 AUC, B1 AUC, sample size) corresponding to the integral Betti signatures (the Area Under the Curve of the $0$-th and $1$-st Betti curves), along with the dimension of the manifold (sample size). We plot the features corresponding to the three datasets described in Section~\ref{sec:data}, and to the random matrices, geometric (Euclidean, spherical, hyperbolic) matrices, and  correlation matrices. 
All datasets have TS length up to $2^{15}$, giving an overview of the asymptotic behaviour of the associated Betti curves.  The random/geometric curves are very distinct from each other; correlation is close to them for short samples, but converges to random  for longs samples.

\subsubsection*{Random correlation and spherical geometry} We observe that the Betti curves of spherical and correlation matrices are similar. This is because correlation matrices have an intrinsic spherical geometry. Indeed, the Pearson correlation of $n$-dimensional random variables yields the normalized scalar product between the corresponding vectors in~$\R^{n-1}$; from which one obtains the angle between such vectors by application of the function~$\mathrm{arccos}$. Hence, it is equivalent to the spherical distance between them, explaining  the similar behaviour of their Betti curves shown in Figures~\ref{fig:Drawing} and~\ref{fig:VDL2d_REALvsSim}.

\subsubsection*{Hyperbolic geometry}
We complete our analysis of Betti curves of geometric matrices by investigating also the behaviour of point clouds distributed in the hyperbolic space~$\bH^n$. The choice of radii was $R=0.01,0.05,0.1,0.5,0.7,1,10$. 
On the technical note, the random points in the hyperbolic space $\bH^n$ have been generated in the 
Poincar\'e disc model. For a given $R\in (0,\infty)$ we have randomly generated vectors in the ball $B^n_R$ of vectors of norm $\leq R$. 
The obtained points belong to the hyperbolic space of radius $R$; we have then projected such points on the standard hyperbolic space $\bH^n$ by applying the transformation 
$
r\mapsto (\cosh{\rho(r)} - 1)/(2 + \cosh{\rho(r)})
$; see also the Appendix.

For such configurations, for small radii we observe that the Betti numbers~$\beta_p$, for $p>0$, are close to $0$ (see Figure~\ref{fig:Drawing}). When the radius increases, also B0AUC and B1AUC increase, showing characteristics closer to random correlation. This effect holds consistently across the choice of dimension, albeit more visible for higher dimensions (see~Figure~\ref{fig:VDL2d_REALvsSim}D).   
Note also from Figure~\ref{fig:Drawing} that different radii have different behaviours of $0$-th Betti curves, with lower slope for smaller radii.


\subsection{Curvature of real-world data}

\begin{figure*}[h!]
\begin{center}
        \includegraphics[width=1.0\textwidth]{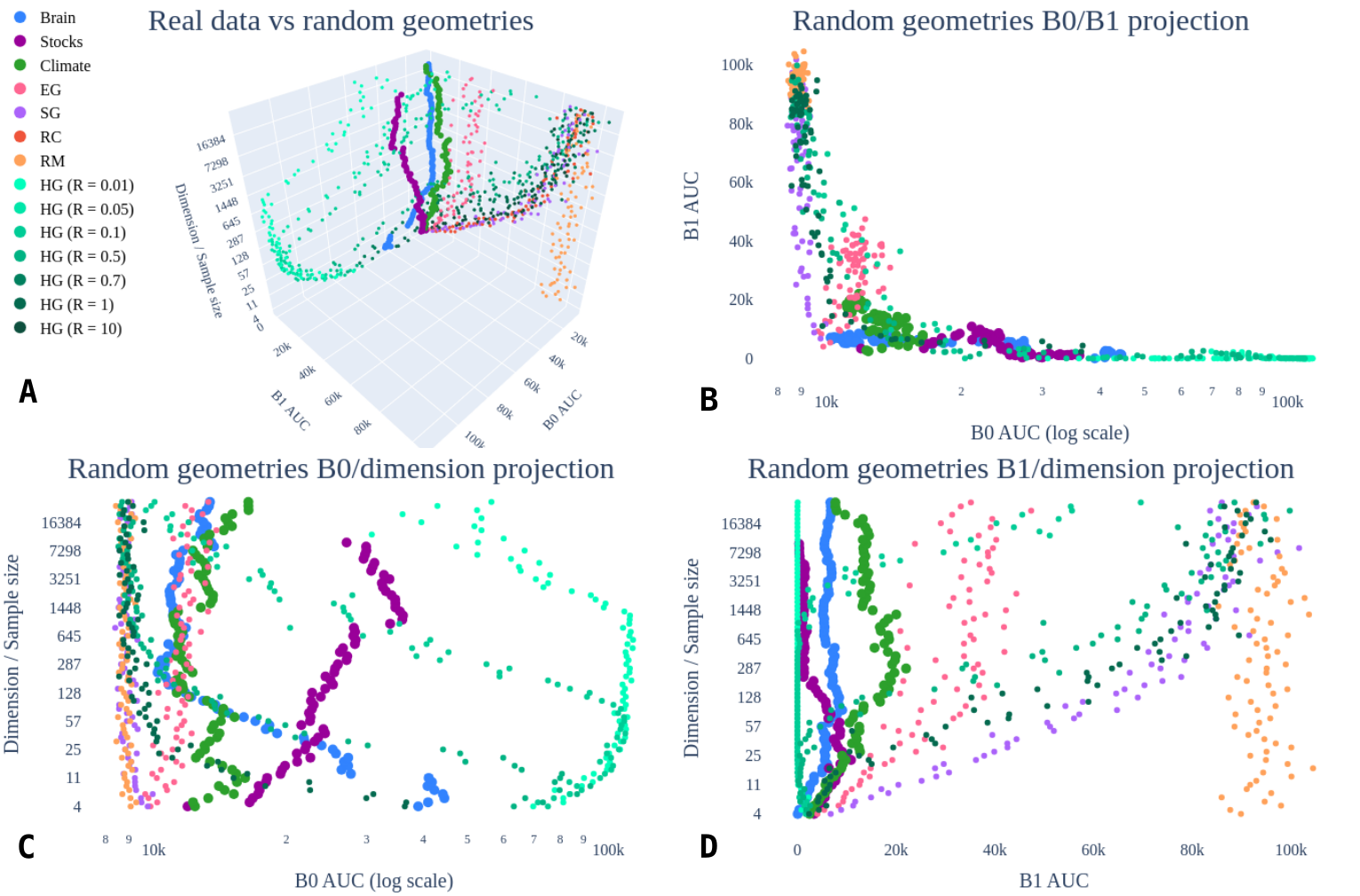}        
	\caption{ 
 Integral Betti signatures 
 for real data correlation matrices, random matrices (RM)/Euclidean geometry~(EG)/spherical geometry (SG)/hyperbolic geometry (HG)/random correlation (RC) matrices. Subfigure A shows the 3D plot of integral Betti signatures for different sample sizes/dimensionality for the real data and the simulated random geometries. Subfigures B, C and D present the 2D projections B0AUC/B1AUC, B0AUC/dimensionality and B1AUC/dimensionality, respectively.}
	\label{fig:VDL2d_REALvsSim}
\end{center}
\end{figure*}

We analyse the case of three main sources of different nature: brain, stocks, and climate data; these are described in Section~\ref{sec:data}. In order to analyse the influence of the sample size on the associated Betti curves, for every dataset, we start with an initial segment of length $4$ (extremely short and not meaningful in practice, but illustrative concerning the general effects), and then take longer segments until $2^{15}$ points (very long, longer than used in many practical applications). We also consider (simulated) random and geometric matrices as comparison models; the size of these matrices is also fixed to $90\times 90$, with sample points taken from manifolds of dimensions corresponding to the sample size.
 In Figure~\ref{fig:VDL2d_REALvsSim} we show the associated 3-dimensional plot. 

 As discussed in the previous subsections, the Betti curves of random matrices are quite distinct from those of any of the geometric matrices, at least for low dimensions (sample size). 
 For small radii, the hyperbolic matrices are close to 0 in the B1AUC-coordinate -- related to the fact that for small radii, there are few $1$-dimensional cycles. On the other hand, for large radii the Betti curves of the hyperbolic geometric matrices are comparable with those of random correlation and spherical geometric matrices -- particularly so in high dimensions. 
 
 For random matrices, the point cloud is approximately normally distributed in the B0AUC-B1AUC projection, and there is no substantial dependence on the dimensionality. In contrast, the Euclidean geometry shows a 
 less pronounced but still present systematic effect of sample size. 

Turning to real data, we observe that its topological features sit clearly between hyperbolic (of small radii) and Euclidean features. The longer the time series, the more stable the Integral Betti signatures 
become, and has a smaller value of B1AUC (fewer $1$-dimensional cycles). 
Particularly, the stock data show few $1$-dimensional cycles represented in the B1AUC coordinates, and it seems to be closer to the hyperbolic distributions (with small radii). Also, the behaviour in B0AUC changes with the sample size. 
Note that the behaviour of the stock data shown in Figure~\ref{fig:VDL2d_REALvsSim}  agrees with previous discussions on the hyperbolicity of the stock market~\cite{hypfinance} 
and that this is also consistent with part of the results of~\cite{TDA_geometry_sptrain}.
We point out that both the stock market ~\cite{Hartman2018} and climate~\cite{Hlinka2014} time series contain non-negligible nonstationarities, while fMRI data are concatenated across subjects; this may contribute to some visible jumps in the integral Betti signatures when changing time series length. This gets smoothed out when a stationary model of the stock data is sampled (simulations not shown), making even clearer the general drift with increasing sample size towards an apparent low-radius hyperbolic geometry.

Although possibly related to an intrinsic hyperbolicity of the system, our results might also be related to other characteristics of the datasets; e.g.~to  a small rank property, whose effect on the Betti curves was theoretically investigated in~\cite{curto2021betti}.  In the last part of our experiments, we focus on this effect.


\subsection{Modularity}
\label{sec:modularity}

\begin{figure}[h]
\begin{center}
    \includegraphics[width=\columnwidth]{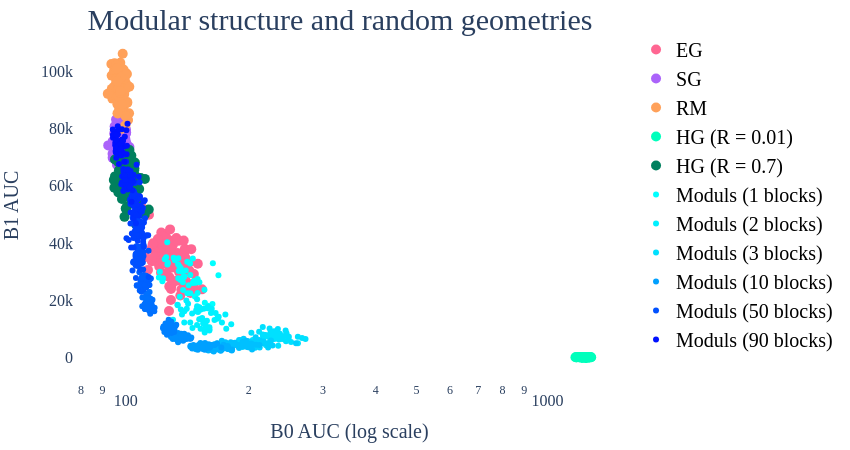}
	\caption{Integral Betti signatures 
 of geometric distributions and modular systems (for a growing number of modules). Shortcuts: EG - Euclidean geometry, SG - spherical geometry, RM - random matrices, HG - hyperbolic geometry. The sample size/dimensionality is 400. }
	\label{fig:ModularSystems}
\end{center}
\end{figure}

We aim   to understand the dependence of the Betti curves with respect to (small) rank matrices.
In order to create (small) rank matrices we have constructed modular systems with modules of the same size. 
To do so, we start with $90$ independent random  time series (white noise) with the length of $400$ sample points; for constructing a modular system of size $m\leq 90$, we first  select $m$ of such time series; say $M_1,\dots,M_m$. We consider the vector $[M_1,\dots,M_m,M_1,\dots,M_m,\dots,M_1,\dots,M_{k}]$ where $k \equiv 90 \mod m$, and in total we create $90$ time series with $m$ modules of approximately the same size. After copying the time series, we add a small-amplitude white noise. The effect of this modular structure -- projected onto the B0AUC-B1AUC plane 
 -- is shown in Figure~\ref{fig:ModularSystems}; here, we plot features for the random and geometric distributions, along with those associated to the described matrices with $m$ modules, for $m=1,\dots , 90$. 
The modules show a gradient starting close to Euclidean (corresponding to 1 module) and then moving to hyperbolic (until about 5 modules). The gradient stays in the ``hyperbolic regime'' (at different radii) and then continues towards random (when reaching the full rank).  

This simulation demonstrates that, depending on the number of modules,  the Betti curves in the (B0AUC, B1AUC) coordinates yield features which are similar to certain specific geometric matrices; going from Euclidean to hyperbolic of different radii to spherical to random.
When the modular system is described by a single module, this behaviour is explained by the fact that, locally, the spherical geometry looks like Euclidean in one less dimension (consider points sampled close to a pole of a sphere effectively close to Gaussian sample from the tangential plane). When the number of modules increases (but still far from yielding a full rank matrix), the shift to the hyperbolic direction of small radii (hence, of few 1-cycles) can be explained by the features being ``governed'' by the behaviour of only a small amount of points on the sphere, the centroids provided by the modules themselves; there is a shift to bigger B0AUC, determined by the fact that, overall, the number of sample points is still higher than the number of centroids. When the number of modules goes to the full matrix rank, the distribution on the sphere tends to be more homogeneous, hence behaving like the random distribution. This effect depends on the dimension of the space.


\section{Discussion}

Betti curves  are  well understood from a theoretical point of view, and were shown to provide an effective tool in distinguishing random  from (Euclidean) geometric matrices~\cite{Giusti13455}. The asymptotic behaviour of the Betti curves associated to long time series was also  investigated~\cite{Giusti13455}. However, in practical applications, one commonly deals with short time series. We showed that the qualitative behaviour of Betti curves may change considerably.  
We also extended  the   analysis previously developed in~\cite{Giusti13455} (in the case of random and Euclidean geometry) to spherical and hyperbolic geometries. We have used the AUC in conjunction with the $0$ and $1$-dimensional Betti curves as main qualitative measures, which we referred to as \emph{integral Betti signatures}. The  outputs depend on several parameters; first, the number of points sampled and the dimension of the ambient manifolds. We investigated the dependence on these parameters, showing that the Betti curves distinguish the three geometric models. We have considered the three geometries as base models in the subsequent analysis of real data.

Applying Pearson correlation, we  observed that Betti curves associated to real world time series 
  sit in between hyperbolic far from being spherical. On the other hand, for larger radii, hyperbolic geometric matrices seem to behave like spherical ones. This is intriguing, given that truly random correlation matrices behave rather like spherical geometric matrices.   Data indeed can have true hyperbolic underlying geometry, as shown across contexts~\cite{hypfinance,doi:10.1126/sciadv.aaq1458,10.3389/fphy.2018.00007}, but this ``hyperbologenic effect'' can also be an artefact due to small data rank. Following the theoretical study in~\cite{curto2021betti}, we analysed the small rank effect on Betti curves in simulations, showing how the qualitative trends depend on the number of modules. 
Moreover, we observed the role of the choice of the metric adopted in the construction of the order complexes. In Appendix Figure~C.5, we show how the Betti curves of the random geometric matrices change when using a different metric (defined on the same underlying manifold). Also, the ordering of the filtration in the construction of the order complexes matters. This shows again how the metric (rather than solely the distribution of sample points) plays an important role in the computation of Betti curves. For real data (see Figure~C.6), the dependence on the distribution of points is deeper -- e.g. for climate and brain. The main differences are visible for \emph{shorter} time series, for lengths typical for real-world scenarios. Finally, we study the Betti numbers given by the Pearson correlation as a measure of dependence, a choice shown to capture sufficiently pairwise dependence in  brain~\cite{Hlinka2011}, climate~\cite{Hlinka2014} and financial~\cite{Hartman2018} data. While the correlation matrix captures rich (multivariate) global structure, high-order interactions are of interest~\cite{Rosas2022}, although e.g.~in the case of brain not dominant~\cite{Martin2017}. 

We now turn to our main motivation. When dealing with complex networks, it is believed that most complex networks have some hidden network geometry~\cite{bianconi}. 
There is a general belief that further investigations of such hidden geometries will deepen our understanding of the fundamental laws describing relationships between structure and function of complex networks~\cite{navigability}. However, connections between the hidden geometry and the combinatorial properties of the networks are not yet totally understood. Some works have  partially described  these relations, e.g.~\cite{bianconi} in the context of emergent hyperbolic geometry in growing simplicial complexes. However, a complete picture (for real data) is still missing. 
Persistent homology and Betti curves have proved to be great qualitative tools in investigating related questions. 
In fact, geometric properties, like the curvature, have effects on the homology of the underlying manifolds~\cite{goldberg1998curvature}, hence, on persistent homology. Moreover, PH-derived features  are easily computable, thanks to various algorithms freely available. However, more relations between persistent homology and the curvature of data manifolds have yet to be discovered; see the recent advances~\cite{ curto2021betti,Lim2020,Bubenik_2020, otter_curv, hacquard2023euler}.
 Our analysis is meant to contribute to the advancing of our understanding of the hidden structure of data manifolds; we believe that PH- and TDA-based approaches will be beneficial for these tasks, as we qualitatively showed in our results. 

\subsection*{Conclusions}

As highlighted earlier, high-dimensional datasets are important in many scientific domains, including biology, economics or climate. Data 
are increasingly assumed to lie on  hidden geometric manifolds. Therefore, understanding geometric and topological properties, such as the curvature, of these manifolds can be beneficial for uncovering hidden structures of data.  Estimations of the 
curvature have been previously considered with analytical methods. In this work, we further developed the integral Betti signature approach based on modern tools from Algebraic Topology and Topological Data Analysis, documenting its ability to characterize the differences, as well as interesting similarities, among topologies of a range of random and geometric distributions. 
Importantly, we showed that these relations depend crucially and systematically on further aspects including sample size/space dimension, hyperbolic space radius, and artificially imposed modularity; features that need thus to be taken into account when interpreting observed Betti signatures.
Finally, application to real-world data demonstrated that
brain, climate, and financial datasets tend to be the closest to hyperbolic, i.e.~fast-transmission, networks.

\subsection*{ ACKNOWLEDGMENTS.}
The authors acknowledge the support of the Czech Science Foundation grant 21-17211S, JH acknowledges the support 
of the CSF grant 21-32608S, AP and JH acknowledge the support of the Ministry of Health Czech Republic – DRO 2021 (``NIMH, IN: 00023752''). The publication was also supported by ERDF-Project Brain dynamics, No. CZ.02.01.01/00/22\_008/0004643. 

\bibliographystyle{abbrv}
\bibliography{biblio}

\def\cprime{$'$} \providecommand{\MR}[1]{}\providecommand{\MRhref}[1]{}
\begin{thebibliography}{10}

\bibitem{YahooFinance2023}
{Yahoo! Finance}.
\newblock https://finance.yahoo.com/, 2023.
\newblock Accessed: 2023-10-06.

\bibitem{AlanisLobato2016DistanceDB}
G.~Alanis-Lobato and M.~Andrade.
\newblock Distance distribution between complex network nodes in hyperbolic
  space.
\newblock {\em Complex Syst.}, 25, 2016.

\bibitem{bianconi}
G.~Bianconi and C.~Rahmede.
\newblock Emergent hyperbolic network geometry.
\newblock {\em Sci. Rep.}, 7, 2017.

\bibitem{navigability}
M.~Boguñá, D.~Krioukov, and K.~Claffy.
\newblock Navigability of complex networks.
\newblock {\em Nature Phys}, 5, 2009.

\bibitem{MR1744486}
M.~R. Bridson and A.~Haefliger.
\newblock {\em Metric spaces of non-positive curvature}, volume 319 of {\em
  Grundlehren der mathematischen Wissenschaften}.
\newblock Springer-Verlag, Berlin, 1999.

\bibitem{Bubenik_2020}
P.~Bubenik, M.~Hull, D.~Patel, and B.~Whittle.
\newblock Persistent homology detects curvature.
\newblock {\em Inverse Problems}, 36(2):025008, 2020.

\bibitem{ourpaper}
L.~Caputi, A.~Pidnebesna, and J.~Hlinka.
\newblock Promises and pitfalls of topological data analysis for brain
  connectivity analysis.
\newblock {\em NeuroImage}, 238:118245, 2021.

\bibitem{MR3715451}
G.~Carlsson, T.~Ishkhanov, V.~de~Silva, and A.~Zomorodian.
\newblock On the local behavior of spaces of natural images.
\newblock {\em Int. J. Comput. Vis.}, 76(1):1--12, 2008.

\bibitem{curto2021betti}
C.~Curto, J.~Paik, and I.~Rivin.
\newblock Betti curves of rank one symmetric matrices.
\newblock In F.~Nielsen and F.~Barbaresco, editors, {\em Geometric Science of
  Information}, pages 645--655, 2021.

\bibitem{vmv.20171272}
T.~K. Dey, S.~Mandal, and W.~Varcho.
\newblock {Improved Image Classification using Topological Persistence}.
\newblock In M.~Hullin, R.~Klein, T.~Schultz, and A.~Yao, editors, {\em Vision,
  Modeling and Visualization}, 2017.

\bibitem{MR3522608}
C.~Fefferman, S.~Mitter, and H.~Narayanan.
\newblock Testing the manifold hypothesis.
\newblock {\em J. Amer. Math. Soc.}, 29(4):983--1049, 2016.

\bibitem{fernandez2023intrinsic}
X.~Fern{\'a}ndez, E.~Borghini, G.~Mindlin, and P.~Groisman.
\newblock Intrinsic persistent homology via density-based metric learning.
\newblock {\em Journal of Machine Learning Research}, 24(75):1--42, 2023.

\bibitem{gidea}
M.~Gidea.
\newblock Topology data analysis of critical transitions in financial networks.
\newblock {\em SSRN Electronic Journal}, 01 2017.

\bibitem{Giusti13455}
C.~Giusti, E.~Pastalkova, C.~Curto, and V.~Itskov.
\newblock Clique topology reveals intrinsic geometric structure in neural
  correlations.
\newblock {\em Proceedings of the National Academy of Sciences},
  112(44):13455--13460, 2015.

\bibitem{goldberg1998curvature}
S.~I. Goldberg.
\newblock {\em Curvature and homology}.
\newblock Courier Corporation, 1998.

\bibitem{TDA_geometry_sptrain}
A.~Guidolin, M.~Desroches, J.~D. Victor, K.~P. Purpura, and R.~S.
\newblock Geometry of spiking patterns in early visual cortex: a topological
  data analytic approach.
\newblock {\em J. R. Soc. Interface}, 19(196):20220677, 2022.

\bibitem{hacquard2023euler}
O.~Hacquard and V.~Lebovici.
\newblock Euler characteristic tools for topological data analysis, 2023.

\bibitem{Hartman2018}
D.~Hartman and J.~Hlinka.
\newblock Nonlinearity in stock networks.
\newblock {\em Chaos: An Interdisciplinary Journal of Nonlinear Science},
  28(8):083127, 2018.

\bibitem{hatcher}
A.~Hatcher.
\newblock {\em {Algebraic topology}}.
\newblock Cambridge Univ. Press, Cambridge, 2000.

\bibitem{Hlinka2014}
J.~Hlinka, D.~Hartman, and M.~et~al. Vejmelka.
\newblock {Non-linear dependence and teleconnections in climate data: sources,
  relevance, nonstationarity}.
\newblock {\em Climate Dynamics}, 42(7-8):1873--1886, 2014.

\bibitem{Hlinka2013}
J.~Hlinka, D.~Hartman, M.~Vejmelka, J.~Runge, N.~Marwan, J.~Kurths, and
  M.~Palu\v{s}.
\newblock Reliability of inference of directed climate networks using
  conditional mutual information.
\newblock {\em Entropy}, 15(6):2023--2045, 2013.

\bibitem{Hlinka2011}
J.~Hlinka, M.~Palu{\v{s}}, M.~Vejmelka, D.~Mantini, and M.~Corbetta.
\newblock {Functional connectivity in resting-state fMRI: Is linear correlation
  sufficient?}
\newblock {\em NeuroImage}, 54(3):2218--2225, 2011.

\bibitem{KAHLE20091658}
M.~Kahle.
\newblock Topology of random clique complexes.
\newblock {\em Discrete Mathematics}, 309(6):1658--1671, 2009.

\bibitem{articleKahle}
M.~Kahle and E.~Meckes.
\newblock Limit theorems for betti numbers of random simplicial complexes.
\newblock {\em Homology, Homotopy and Applications}, 15(1):343--374, 2013.

\bibitem{Kalnay1996}
E.~Kalnay, M.~Kanamitsu, R.~Kistler, W.~Collins, D.~Deaven, L.~Gandin,
  M.~Iredell, S.~Saha, G.~White, J.~Woollen, Y.~Zhu, M.~Chelliah, W.~Ebisuzaki,
  W.~Higgins, J.~Janowiak, K.~C. Mo, C.~Ropelewski, J.~Wang, A.~Leetmaa,
  R.~Reynolds, R.~Jenne, and D.~Joseph.
\newblock {The {NCEP/NCAR} 40-year reanalysis project}.
\newblock {\em Bulletin of the American Meteorological Society},
  77(3):437--471, 1996.

\bibitem{hypfinance}
M.~Keller-Ressel and S.~Nargang.
\newblock {The Hyperbolic Geometry of Financial Networks}.
\newblock {\em Sci Rep}, 11:4732, 2021.

\bibitem{Kopal2020}
J.~Kopal, A.~Pidnebesna, D.~Tome\v{c}ek, J.~Tint\v{e}ra, and J.~Hlinka.
\newblock {Typicality of Functional Connectivity robustly captures motion
  artifacts in rs-fMRI across datasets, atlases and preprocessing pipelines}.
\newblock {\em Human Brain Mapping}, 41(18):5325--5340, 2020.

\bibitem{PhysRevE.82.036106}
D.~Krioukov, F.~Papadopoulos, M.~Kitsak, A.~Vahdat, and M.~Bogu\~n\'a.
\newblock Hyperbolic geometry of complex networks.
\newblock {\em Phys. Rev. E}, 82:036106, Sep 2010.

\bibitem{lee-discriminativePH}
H.~Lee, M.~Chung, H.~Kang, B.-N. Kim, and D.~Lee.
\newblock Discriminative persistent homology of brain networks.
\newblock {\em Proceedings - International Symposium on Biomedical Imaging},
  pages 841--844, 03 2011.

\bibitem{Lim2020}
S.~Lim, F.~Memoli, and O.~B. Okutan.
\newblock {Vietoris-Rips Persistent Homology, Injective Metric Spaces, and The
  Filling Radius}, 2020.

\bibitem{geommmod}
B.~Lin, X.~He, and J.~Ye.
\newblock A geometric viewpoint of manifold learning.
\newblock {\em Appl. Inform.}, 2(3), 2015.

\bibitem{Rmanlearn}
T.~Lin and H.~Zha.
\newblock {Riemannian manifold learning}.
\newblock {\em IEEE Trans. Pattern Anal. Mach. Intell.}, 30(5), 2008.

\bibitem{Martin2017}
E.~Martin, J.~Hlinka, A.~Meinke, F.~D{\v{e}}cht{\v{e}}renko, J.~Tint{\v{e}}ra,
  I.~Oliver, and J.~Davidsen.
\newblock {Network Inference and Maximum Entropy Estimation on Information
  Diagrams}.
\newblock {\em Scientific Reports}, 7(1), 2017.

\bibitem{munkres}
J.~R. Munkres.
\newblock {\em Elements of Algebraic Topology}.
\newblock Addison Wesley Publishing Company, Inc., 2725 Sand Hill Road Menlo
  Park, California 94025, 1984.

\bibitem{Rosas2022}
F.~E. Rosas, P.~A.~M. Mediano, A.~I. Luppi, T.~F. Varley, J.~T. Lizier,
  S.~Stramaglia, H.~J. Jensen, and D.~Marinazzo.
\newblock {Disentangling high-order mechanisms and high-order behaviours in
  complex systems}.
\newblock {\em Nature Physics}, 18(5):476--477, may 2022.

\bibitem{SHARPEE2019101}
T.~O. Sharpee.
\newblock An argument for hyperbolic geometry in neural circuits.
\newblock {\em Current Opinion in Neurobiology}, 58:101--104, 2019.

\bibitem{doi:10.1073/pnas.2100473118}
D.~Sritharan, S.~Wang, and S.~Hormoz.
\newblock {Computing the Riemannian curvature of image patch and single-cell
  RNA sequencing data manifolds using extrinsic differential geometry}.
\newblock {\em Proceedings of the National Academy of Sciences},
  118(29):e2100473118, 2021.

\bibitem{10.1145/3394486.3403224}
P.~Tabaghi and I.~Dokmani\'{c}.
\newblock Hyperbolic distance matrices.
\newblock {\em Proceedings of the 26th ACM SIGKDD International Conference on
  Knowledge Discovery and Data Mining}, page 1728–1738, 2020.

\bibitem{10.3389/fphy.2018.00007}
B.~Tadić, M.~Andjelković, and M.~Šuvakov.
\newblock Origin of hyperbolicity in brain-to-brain coordination networks.
\newblock {\em Frontiers in Physics}, 6, 2018.

\bibitem{otter_curv}
R.~Turkes, G.~F. Montufar, and N.~Otter.
\newblock On the effectiveness of persistent homology.
\newblock In S.~Koyejo, S.~Mohamed, A.~Agarwal, D.~Belgrave, K.~Cho, and A.~Oh,
  editors, {\em Advances in Neural Information Processing Systems}, volume~35,
  pages 35432--35448. Curran Associates, Inc., 2022.

\bibitem{Tzourio-Mazoyer2002}
N.~Tzourio-Mazoyer, B.~Landeau, D.~Papathanassiou, F.~Crivello, O.~Etard,
  N.~Delcroix, B.~Mazoyer, and M.~Joliot.
\newblock {Automated Anatomical Labeling of Activations in SPM Using a
  Macroscopic Anatomical Parcellation of the MNI MRI Single-Subject Brain}.
\newblock {\em NeuroImage}, 15(1):273--289, 2002.

\bibitem{TDAreview1}
L.~Wasserman.
\newblock Topological data analysis.
\newblock {\em Annual Review of Statistics and Its Application}, 5(1):501--532,
  2018.

\bibitem{Whi2021.03.25.436730}
W.~Whi, S.~Ha, H.~Kang, and D.~S. Lee.
\newblock {Hyperbolic disc embedding of functional human brain connectomes
  using resting state fMRI}.
\newblock {\em Netw Neurosci}, 6:745--–764, 2022.

\bibitem{doi:10.1126/sciadv.aaq1458}
Y.~Zhou, B.~H. Smith, and T.~O. Sharpee.
\newblock Hyperbolic geometry of the olfactory space.
\newblock {\em Science Advances}, 4(8), 2018.

\end{thebibliography}
\newpage
\appendix

\onecolumn
\section{Betti curves of  order complexes}\label{sec:Betti curves}

In this section, we briefly recall the topological pipeline introduced by Giust \emph{et al.~}in \cite{Giusti13455}. To each symmetric matrix~$M$, we associate specific topological features called \emph{Betti curves}, one for each dimension $n\in \N$. Roughly, the Betti curve in dimension~$n$ describes the  number of $n$-dimensional holes of a graph with adjacency matrix~$M$, at various thresholds, as a function of the edge density. Following \cite{Giusti13455, curto2021betti}, we now proceed with a detailed description of the pipeline.

\subsection{The order complex}
Recall that a \emph{graph} $G$ is a pair $G=(V,E)$ given by a finite set $V$, whose elements are called vertices, and a set $E\subseteq V\x V$; each element $e$ of $E$, called an edge, is described by an unordered pair  of  
distinct vertices $\{v,w\}$, which are the endpoints of $e$. We observe here that multiple edges between two vertices, and self-loops (i.e.~edges of type $\{v,v\}$) are not allowed. We will only deal with finite graphs, which means, the set of vertices $V$ is finite.  For given graphs $G_1=(V_1,E_1)$ and $G_2=(V_2,E_2)$, we say that $G_1$ is a subgraph of $G_2$, and we write $G_1\subseteq G_2$, if  $V_1$  is a subset of $V_2$, and if $E_1$ is a subset of~$E_2$.

\begin{rem}\label{rem:hfhf}
	Let $G=(V,E)$ be a graph and assume  $V$ to be of cardinality $n$. Consider a bijective function $f\colon \{1,\dots, n\}\to V$ from the set of natural numbers between $1$ and $n$ to~$V$. Then, we can associate to~$G$ a symmetric matrix $A$ called an adjacency matrix of~$G$. For indices~$i,j\in  \{1,\dots, n\}$, let $A(i,j)\coloneqq 1$  if, and only if, $\{f(i),f(j)\}$ is an edge of $G$, and $0$ otherwise. As edges are given by unordered pairs of vertices, we have  $A(i,j)=A(j,i)$ for all $i$ and $j$ in $\{1,\dots, n\}$, and the procedure yields a symmetric matrix. 
\end{rem}  

As recalled in Remark~\ref{rem:hfhf}, to each graph we can associate an adjacency matrix; analogously, to each symmetric matrix with values in $\{0,1\}$, we can associate a graph. 
More generally, we can consider symmetric  real-valued matrices. In such case, to each $N\x N$ symmetric matrix $M$ with distinct non-zero real-valued entries, we associate a whole family of graphs 
\begin{equation}\label{eqApp:ordcompl}
\mathrm{ord}(M)\coloneqq G_0\subseteq G_1\subseteq \dots\subseteq G_k \ ,
\end{equation}
called  the \emph{order complex}  of $M$ \cite[Def.~2, SI]{Giusti13455}. The construction of $\mathrm{ord}(M)$ starts from a totally disconnected graph~$G_0$ on $N$ vertices, and proceeds step by step by adding new edges, as indicated by the entries of $M$. To be more precise, the construction proceeds as follows.
 Let $k$ be the number of non-trivial off-diagonal entries of $M$, counted without repetitions, and let $[a_1,\dots,a_k]$ be the ordered sequence of such (distinct) real values, sorted in a decreasing order. As a first step, let $G_0$ be the graph on $N$ vertices and no edges, with vertices  ordered from $0$ to $N$. Then, inductively construct $G_s$ from $G_{s-1}$ by adding an edge $\{i,j\}$ to $G_{s-1}$ for indices $i, j\in \{0,\dots,N\}$ such that $M(i,j)=a_s$. This iterative construction describes a family of graphs $G_s$, for $s=0,\dots, k$.
  Observe that the number~$k$ is bounded by $N \choose 2$ -- ${N\choose 2}$ being the number of off-diagonal entries of $M$ -- and that each graph $G_{s-1}$ is a subgraph of~$G_s$, making $\mathrm{ord}(M)$ into a sequence of subgraphs of $G_k$. 
  
  In concrete applications, the $N\times N$ symmetric matrix $M$ has often non-zero off-diagonal elements, and $k$ is then equal to~${N\choose 2}$. In the follow-up, the graphs appearing in an order complex will be always indexed by $$\rho= \frac{s}{{N\choose 2}} \ ,$$  the edge density of $G_k$. 
\begin{rem}
The order complex is invariant under monotonic transformations.  This property allows great flexibility in the applications and, especially in presence of non-linearity, it can be used to detect geometric signatures of structure and/or randomness \cite{Giusti13455}.  Furthermore, the construction depends only on the (combinatorics of the)  symmetric matrix. 
\end{rem}

In the construction, we have sorted the values of $M$ in a decreasing order. Analogously, sorting the values in an increasing order yields another (generally distinct) family of subgraphs of $G_k$. 

\subsection{Clique complexes}
Given a sequence of graphs, it is customary in Topological Data Analysis to construct a sequence of higher dimensional spaces called simplicial complexes. We recall the definition:

 \begin{df}\label{defsimplcompl}
	An \emph{(abstract) simplicial complex} on a set $V$ is a collection $\Sigma$   of  non-empty finite subsets $\sigma$ of  $V$, closed under taking subsets:
if $\sigma\in \Sigma$  and $\tau\subseteq\sigma$ is non-empty, then $\tau\in \Sigma$.
\end{df}

The elements $\sigma$ of $\Sigma$ are called \emph{simplices} and the elements of $V$ also called vertices of $\Sigma$. The dimension $\mathrm{dim}(\sigma)$ of a simplex~$\sigma$ is given by the number of its vertices: $\mathrm{dim}(\sigma)\coloneqq |\sigma|-1$, where, if $\sigma=[v_1,\dots,v_j]$ then $|\sigma|=j$. For example, a vertex $v\in V$, seen as element $\{v\}$ in $\Sigma$, has dimension $0$. The dimension of a simplicial complex is the maximum dimension across its simplices.
Graphs are straightforward examples of $1$-dimensional simplicial complexes:

\begin{rem}
A graph $G=(V,E)$ is, in particular, a simplicial complex whose vertex set is $V$ and the other simplices are given by the edges of $G$. The dimension of $G$, seen as a simplicial complex, is $1$.
\end{rem}

To a sequence of graphs, it is possible to associate various sequences of simplicial complexes. In this work, we consider the so-called clique complexes. Recall that a $n$-clique, in a graph $G$, is a complete sub-graph on $n$ vertices; in particular, vertices are $1$-cliques and  edges are $2$-cliques.

\begin{df}\label{def:clique complex}
	Let  $G=(V,E)$ be a graph.
 The \emph{clique complex} $\widetilde{G}$ associated to $G$ is the simplicial complex on the set $V$, whose simplices are precisely the cliques of $G$.
\end{df}

By definition, the simplices of $\widetilde{G}$ are given by all the complete subgraphs of $G$. Each $n$-clique corresponds to a $(n-1)$-simplex of $\widetilde{G}$. 
\begin{ex}
If $G$ is the complete graph on $3$ vertices $\{v_0,v_1,v_2\}$, then the simplicial complex~$\widetilde{G}$ consists of the vertices $\{v_0\},\{v_1\},\{v_2\}$, of the edges $\{v_0,v_1\}, \{v_1,v_2\} $ and $\{v_0,v_2\}$, together with the $2$-simplex corresponding to the whole clique $\{v_0,v_1,v_2\}$. 
\end{ex}

Every simplicial complex can be realized geometrically -- cf.~\cite{munkres}. In the previous example, this geometric realization can be illustrated as follows:
\begin{center}
	\begin{tikzpicture}[node distance={12mm}, main/.style = {}] 
		\node[main] (0) {$v_0$}; 
		\node[main] (2) [right of=0] {$v_2$}; 
		\node[main] (1) [above of=0] {$v_1$}; 
		\draw (0) -- (1);
		\draw (1) -- (2);
		\draw (0) -- (2);
	\end{tikzpicture} 
	\begin{tikzpicture}
		\tikzstyle{point}=[circle,thick,draw=black,fill=black,inner sep=0pt,minimum width=2pt,minimum height=2pt]
		\node[point] at (1.5,0.7) {};
		\node[point] at (1.5,-0.7) {};
		\node[point] at (3.0,-0.7) {};
		\draw[fill=green,opacity=0.4] (1.5,0.7) -- (1.5,-0.7) -- (3.0,-0.7) --cycle; 
		
		\draw[->] (-0.5,-0.2) -- (1,-0.2);
		
	\end{tikzpicture}
\end{center}
To be more precise, an $n$-simplex $\sigma=\{v_0,\dots,v_n\}$  is geometrically realized as the convex hull of $n+1$ geometrically independent vectors $v_0,\dots, v_n$ in $\R^{n+1}$. A $0$-simplex is a point,  a $1$-simplex is depicted as a segment, a $2$-simplex is  a triangle, and so on. Then, simplices are glued together along common faces.

\begin{ex}\label{ex:graph with hole}
	Consider the graph $G$ on seven vertices as described by the following picture:
	\begin{center}
		\begin{tikzpicture}[node distance={12mm}, main/.style = {}] 
		\node[main] (0) {$v_0$}; 
		\node[main] (2) [right of=0] {$v_2$}; 
		\node[main] (1) [above of=0] {$v_1$}; 
		\node[main] (4) [above of=2] {$v_3$}; 
		\node[main] (5) [right of=2] {$v_4$}; 
		\node[main] (6) [right of=4] {$v_5$}; 
		\node[main] (7) [right of=6] {$v_6$}; 
		\draw (0) -- (1);
		\draw (1) -- (2);
		\draw (0) -- (2);
		\draw (1) -- (4);
		\draw (2) -- (4);
		\draw (2) -- (5);
		\draw (4) -- (6);
		\draw (6) -- (7);
		\draw (5) -- (6);
		\draw (5) -- (7);
	\end{tikzpicture} 
\end{center}
Note that, besides $1$- and $2$-cliques, the graph has also the $3$-cliques $\{v_0,v_1,v_2\}$, $\{v_3,v_1,v_2\}$ and also $\{v_4,v_5,v_6\}$. Therefore, the associated clique complex has three $2$-simplices, and these are glued together along the faces. For example, the simplices corresponding to the cliques $\{v_0,v_1,v_2\}$ and $\{v_3,v_1,v_2\}$ are glued together along the segment $\{v_1,v_2\}$. The geometric realization of the graph $G$ is the following:
\begin{center}
	\begin{tikzpicture}
		\tikzstyle{point}=[circle,thick,draw=black,fill=black,inner sep=0pt,minimum width=2pt,minimum height=2pt]
		\node[point] at (1.5,0.7) {};
		\node[point] at (1.5,-0.7) {};
		\node[point] at (3.0,-0.7) {};
		\draw[fill=green,opacity=0.4] (1.5,0.7) -- (1.5,-0.7) -- (3.0,-0.7) --cycle; 
		\node[point] at (3.0,0.7) {};
		\draw[fill=green,opacity=0.4] (1.5,0.7) -- (3.0,-0.7) -- (3.0,0.7) --cycle; 
		\node[point] at (4.5,-0.7) {};
		\node[point] at (4.5,0.7) {};
		\node[point] at (6.0,0.7) {};
		\draw[fill=green,opacity=0.4] (4.5,-0.7) -- (4.5,0.7) -- (6.0,0.7) --cycle;
		
		\draw (4.5,0.7) -- (3,0.7);
		\draw (4.5,-0.7) -- (3,-0.7);		
	\end{tikzpicture}
\end{center} 
\end{ex} 

\subsection{Homology of complexes}

Classical topological invariants of simplicial complexes are the so-called \emph{homology groups}~\cite{hatcher, munkres}. For a simplicial complex $\Sigma$, 
the $k$-th {homology group} $\mathrm{H}_k(\Sigma)$ of $\Sigma$ 
can be thought of as the set of $k$-dimensional holes of the geometric realization of $\Sigma$. To be more precise, consider  the field $\Z_2$ with two elements (i.e.~$0$ and $1$, with sum and product inherited from the usual sum and product of real numbers, reduced mod $2$) and let $C_p(\Sigma)$ be the free $\mathbb{Z}_2$-vector space whose basis consists of the set of $p$-simplices of  $\Sigma$. For every $p\geq 1$ we define the map 
\[
\partial_p\colon C_p(\Sigma)\to C_{p-1}(\Sigma), \quad \partial_p(\sigma)\coloneqq \sum_{\tau \subseteq \sigma, \tau \in C_{p-1}(\Sigma)}\tau
\]
by sending a $p$-simplex $\sigma$ to a formal sum of all its $(p-1)$-faces.
The map $\partial_0$ is defined as the zero map.  An easy computation shows that the composition $\partial_{p}\circ \partial_{p+1}=0$ is the zero map -- \emph{cf.}~\cite{hatcher} -- hence, the image of $\partial_{p+1}$ in $C_p(\Sigma)$ is contained in the kernel of $\partial_{p}$. The homology groups of the simplicial complex $\Sigma$, with $\mathbb{Z}_2$-coefficients are defined as follows:

\begin{df}\label{defhom}
	Let $\Sigma $ be a simplicial complex and $k\geq 0$ a natural number. The $k$-th \emph{homology group} $\mathrm{H}_k(\Sigma)$ of $\Sigma$
	\[
	\mathrm{H}_k(\Sigma)\coloneqq \ker(\partial_k)/\mathrm{Im}(\partial_{k+1})
	\]
	is defined as the quotient of the kernel of the map $\partial_k$ with the image of $\partial_{k+1}$.
\end{df}
 We refrain here from giving the definition of homology groups with more general coefficients, referring to more classical texts on the subject as \cite{hatcher, munkres}, or also the appendices of \cite{Giusti13455,ourpaper}.

\begin{df}
	Let $\Sigma$ be a simplicial complex. The  dimension of $\mathrm{H}_i(\Sigma)$ (over $\mathbb{Z}_2$) is called the $i$-th \emph{Betti number}  (over $\mathbb{Z}_2$) of~$\Sigma$; the $i$-th Betti number is denoted by $\beta_i(\Sigma)$, or simply by $\beta_i$ when the simplicial complex is clear from the context. 
\end{df}

The group $\mathrm{H}_k(\Sigma)$ is in fact a $\mathbb{Z}_2$-vector space and its dimension over $\mathbb{Z}_2$ is well-defined. Betti numbers describe topological and geometric features of simplicial complexes.
In fact, the number of $k$-holes of (the geometric realization of) a simplicial complex $\Sigma$ corresponds to  the $k$-th  {Betti number} $\beta_k(\Sigma)$ associated this way to $\Sigma$. The $0$-th Betti number $\beta_0$ gives the number of connected components, the $1$-st Betti number $\beta_1$ the number of independent loops and the $2$-nd Betti number $\beta_2$ gives the number of $2$-dimensional spheres embedded in $\Sigma$. 
\begin{ex}
	If $X$ is a point, then its Betti numbers $\beta_n(X)$ are  $0$ for every $n>0$, except for  $\beta_0(X)$ which is  $1$. The $0$-th Betti number of two points is $2$ (corresponding to having $2$ connected components) and $0$ in higher dimensions. If we consider $X$ to be the (geometric realization of the) simplicial complex of Example~\ref{ex:graph with hole}, then the Betti numbers $\beta_n(X)$ are  $0$ for every $n>1$, except for  $\beta_0(X)$ and $\beta_1(X)$, which are both  $1$ -- corresponding to $X$ having a single connected component and a cycle.
\end{ex}

\subsection{Betti curves of order complexes}

We now go back to the sequence of graphs appearing in the order complex associated to a symmetric matrix.
For a given order complex $\mathrm{ord}(M)$ represented as the sequence of graphs in \eqref{eqApp:ordcompl}, the pipeline explained in the previous sections yields a sequence of clique complexes
\begin{equation}\label{eq:cliquecompl}
	\widetilde{G_0}\subseteq \widetilde{G_1}\subseteq \dots\subseteq \widetilde{G_k}
\end{equation}
where also the containments 	$\widetilde{G}_{s-1}\subseteq \widetilde{G}_s$ are preserved.
Furthermore, for each simplicial complex $\widetilde{G_j}$ in the sequence of \eqref{eq:cliquecompl}, we can compute the Betti numbers~$\beta_i(\widetilde{G_j})$, hence we can consider the sequence of $i$-th Betti numbers
  \begin{equation}\label{eqApp:Betti curves}
  	\beta_i(\widetilde{G_0}), \beta_i(\widetilde{G_1}), \dots, \beta_i(\widetilde{G_k})
  \end{equation}
and, analogously, when the sequence is indexed on the edge density. 
\begin{df}
	For a symmetric real-valued matrix $M$, the sequences of Betti numbers described in~\eqref{eqApp:Betti curves} are called the ($i$-th) \emph{Betti curves} of $M$. The index $i$ will be called the \emph{Betti dimension}.
\end{df}
 The Betti curves roughly describe the topological dynamics behind the matrix $M$ and provide new invariants that depend only upon the relative order of the entries of the matrix~\cite{Giusti13455}.

\section{Random, geometric and correlation matrices}

We call \emph{random} any symmetric matrix with identically independent real-valued entries, \emph{geometric} any symmetric matrix which is obtained as the distance matrix of sample points uniformly, randomly, distributed on a manifold, and \emph{correlation matrix} any symmetric matrix obtained as the (Pearson) correlation matrix of given time series.

The Betti curves  can reliably detect both the Euclidean geometry and the random structure of symmetric matrices -- see~\cite{Giusti13455}; in our work we complement this observation by exploring the case of spherical and hyperbolic geometries, together with the case of correlation matrices. We review here their definitions.

\subsection{Random matrices}

By random matrices, we mean  \emph{random symmetric matrices with i.i.d.~entries in $[0,1]$.} Note that these matrices are adjacency matrices of Erd\"{o}s–Rényi graphs. For an  Erd\"{o}s–Rényi graph, or the adjacency matrix of an Erd\"{o}s–Rényi graph, we call \emph{random clique complex} the associated clique complex -- cf.~Definition~\ref{def:clique complex}.

The asymptotic  behaviour of the Betti curves arising from random matrices  has been  theoretically studied, and it is now well understood: for each $k\in \N$, there is an interval, depending on $k$, where the $k$-th Betti numbers of a random clique complex~$X$ are non-zero asymptotically almost surely \cite{KAHLE20091658}. Furthermore,  the expected number $\mathbb{E}[\beta_k(X)]$ has been computed, and it has been shown that $\beta_k(X)$ satisfies a Central Limit Theorem  \cite{articleKahle}. 
%
The theoretical predicted results are asymptotic, and depend on the sizes of the input matrices as well. When the size is small (ca.~$20$) in fact,  the behaviour of the Betti curves changes substantially~\cite{Giusti13455}. 



\subsection{Geometric matrices}

We consider the  metric spaces $\R^n$ (the Euclidean space), $\mathbb{S}^n$ (the sphere), and $\bH^n$ (the hyperbolic space). These are the geometric models of the simply connected $n$-manifolds of constant curvature $0$, $1$ and $-1$, respectively. 
Generally speaking, the geometric properties of a manifold, endowed with different underlying metrics, can be very different. For example,  geodesic triangles in the three models, have  different behaviours -- \emph{cf}.~Figure~\ref{fig:triangles}.


\subsubsection{Geometric (Euclidean) matrices}

Euclidean geometry is the default choice in geometric and machine learning representations. The underlying manifold is the Euclidean space~$\R^n$, endowed with the standard Euclidean distance between vectors; for $x=(x_1,\dots,x_n)$ and $y=(y_1,\dots,y_n)$ vectors in $\R^n$, the Euclidean distance
\[
d_E(x,y)\coloneqq \sqrt{\sum_{i=0}^n (x_i - y_i)^2}
\]
is the classical distance between vectors.

In our experiments, we consider $N$ point samples  uniformly independently identically distributed in $[0,1]^d\subseteq \R^d$. For each such point, we get a vector $x_i\in\R^d$. Then, we compute the symmetric matrix~$E$ consisting of all the Eucledian distances: $E(i,j)=E(j,i)\coloneqq d_E(x_i,x_j)$. As the matrix $E$ is symmetric, we can apply the pipeline, the associated rank matrices, and the Betti numbers by density. 


\subsubsection{Geometric (Spherical) matrices}

In the case of spherical geometry, with curvature $+1$, we consider  points sampled on a sphere~$\mathbb{S}^{d-1}\subseteq\R^{d}$, endowed with the spherical distance
\[
d_S(x,y)\coloneqq \arccos{\langle x,y\rangle}
\]
where $\langle x,y\rangle\coloneqq \sum_{i=1}^{d} x_i y_i$ denotes the scalar product of $x$ and $y$ in $\R^{d}$. Note that this is the distance of the shortest path on the sphere (a geodesic) between $x$ and $y$. 

In our experiments, using a Gaussian normal distribution ($\sim \mathcal{N}(0,1)$), we first sample $N$ random points in $\R^d$ and reduce them to unit vectors -- modulo their norm. We point out here that we cannot use the uniform distribution on $[0,1]^d$ because not spherically symmetric. This process gives $N$ unit vectors lying on the sphere~$\mathbb{S}^{d-1}\subseteq\R^{d}$, hence we   compute the associated matrix $S$ consisting of the relative spherical distances. We finally take the associated rank matrix and compute the average Betti curves.





\subsubsection{Geometric (Hyperbolic) matrices }

Hyperbolic geometry, is a non-Euclidean geometry, in which geodesics tend to diverge (as opposed to the spherical one where geodesics tend to converge). There are many models for hyperbolic spaces with negative constant curvature. In the following, we consider the Poincar\'{e} disc model. 
The underlying set of the $d$-dimensional Poincar\'{e} disc model $\bH^d$ is the standard open $d$-dimensional ball 
\[
B^d\coloneqq \{ x=(x_1,\dots,x_d)\in\R^d : |x|=\sum_{i=1}^d x_i< 1\} 
\]
in $\R^d$, endowed with the distance 
\[
d_{\bH^n}(v,w)\coloneqq \arccosh\left( 1 + 2 \frac{\lVert v-w \rVert^2 }{(1- \lVert v\rVert^2)(1- \lVert w\rVert^2)} \right) \ ,
\]
where $\lVert x \rVert$ denotes the standard Euclidean norm of a vector $x$ in $\R^d$. We refer to  \cite{10.1145/3394486.3403224} for an overview of this and other models for the hyperbolic space. 

In order to get random distributions of points in the hyperbolic space we use the approximation of the distribution from~\cite{AlanisLobato2016DistanceDB} at given radii. 
The random sample points have been generated in the $n$-dimensional ball $B^n$ (of vectors of norm $\leq 1$) -- which means, in the Poincar\'e disc model. For a given $R\in (0,\infty)$ we have randomly generated vectors in the ball $B^n_R$ of vectors of norm $\leq R$ by first using a uniform distribution on the sphere~$S^{n-1}$ seen as the boundary of the $n$-dimensional ball; then, we have used the distribution
\[
\rho(r)=\frac{\sinh^{n-1}{r}}{\cosh^{n-1}{(R-1)}}
\]
for getting random distances from $0$. The obtained points belong to the hyperbolic space of radius $R$; we have then projected such points on the standard hyperbolic space $\bH^n$ by applying the transformation 
\[
r\mapsto (\cosh{\rho(r)} - 1)/(2 + \cosh{\rho(r)})\ . 
\]
Using these sample points, we have then computed the hyperbolic distances between them. These matrices are  symmetric and real-valued.


\subsection{Correlation matrices}

The last family of matrices that we investigate is given by correlation matrices. In fact, these are among the main sources of symmetric matrices in applications. Note that correlation matrices can be seen as geometric spherical matrices. We  compute the Betti numbers of correlation matrices (Pearson) on $N$ time points, with uniform distribution in $[-1,1]$. Matrices are then obtained as Pearson correlation of the time series. 




\section{Additional results and notes: different distance measures}

\begin{figure*}[h!]
\begin{center}
		\includegraphics[width=0.49\textwidth]{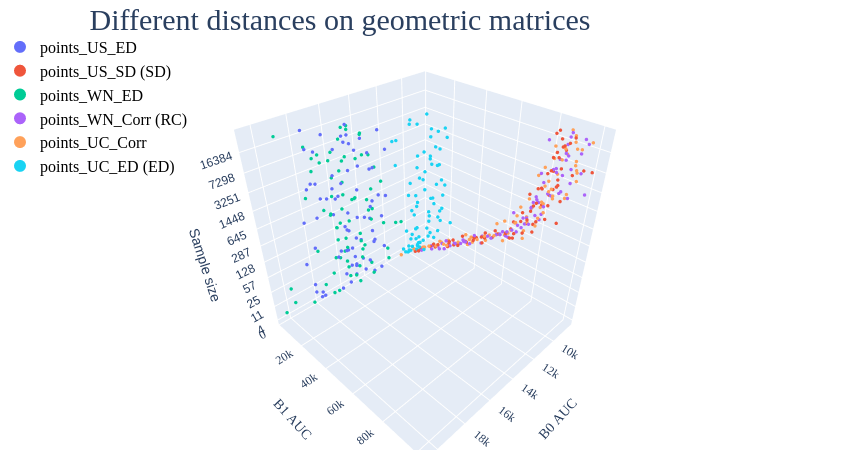}
        \includegraphics[width=0.49\textwidth]{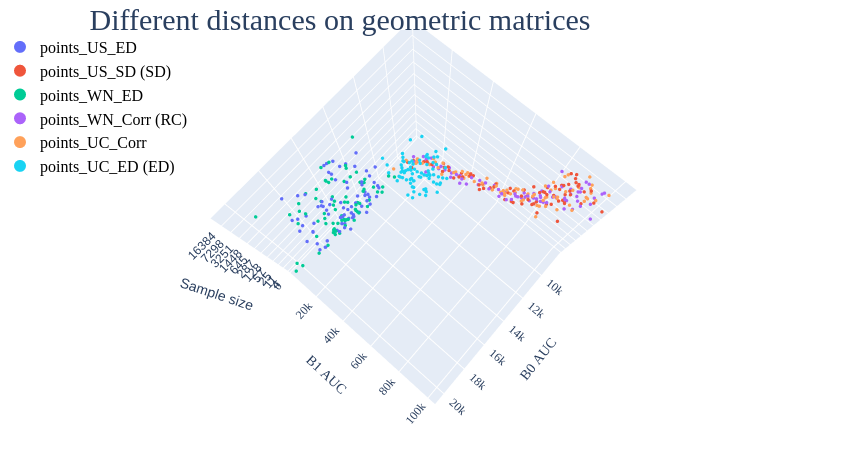}		\label{fig:supp_distanceMeasures}
	\caption{Comparison of integral Betti signatures, 
 computed for random points in the cube/sphere with different distance measures. ED denotes Euclidean distance; Corr stands for Pearson's correlation; WM means white noise (standard multivariate normal distribution); UC denotes uniformly distributed points in a hyper-cube and US - uniformly distributed points on a sphere.  
	Observe that for correlation-based measure, the distribution of points has no effect in these particular cases. 
	Euclidean distance matrices, computed on a white noise, have similar topology as Euclidean distance matrices, computed from the points, uniformly distributed on a sphere. 
	However, they show more hyperbolic properties compared to the Euclidean distance matrices, computed from the points, uniformly distributed in a hyper-cube.}
	
\end{center}
\end{figure*}


\begin{figure*}[h!]
\begin{center}

		\includegraphics[width=0.49\textwidth]{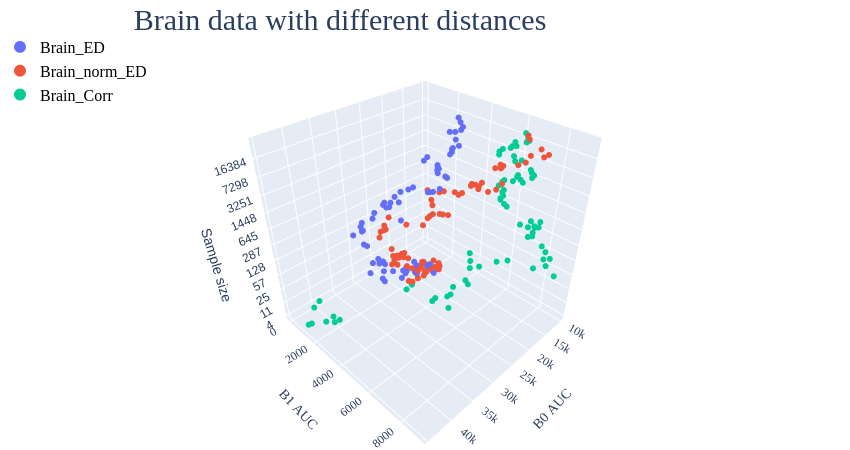}
		\includegraphics[width=0.49\textwidth]{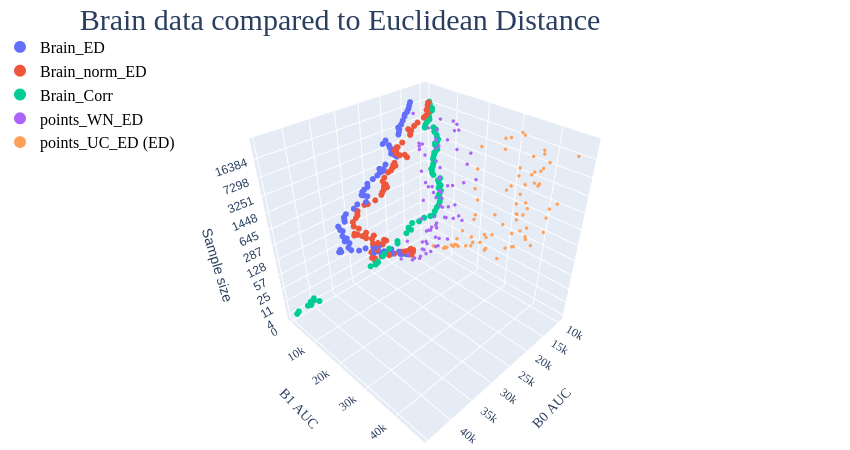}
    \label{fig:supp_distanceMeasures_fMRI}
	\caption{
	Demonstration of the preprocessing influence on integral Betti signatures on the example of fMRI data. Normalization of the time series to variance 1, shifts points from hyperbolic towards euclidean geometry. 
	Distance, computed for the data, is essential, especially for short time series, while long time series results seem to show relatively similar properties. 
	Note also that correlation, computed from brain (fMRI) data, lies close to Euclidean distance, calculated for the white noise, however, is far from Euclidean distance of uniformly distributed points in a hyper-cube and from correlation computed on white noise. This again highlights our observation that similar properties in Betti curves should not be interpreted as necessary similarity in topology. 
    Shortcuts: norm denotes normalization of the time series to variance 1; ED denotes Euclidean distance; Corr stands for Pearson's correlation; WM means white noise (standard multivariate normal distribution); UC denotes uniformly distributed points in a hyper-cube.  	}
\end{center}
\end{figure*}

\end{document}